\documentstyle{amsppt}
\pagewidth{5.05in} \pageheight{8in} \NoRunningHeads \magnification
= 1200

\topmatter
\title Bernstein theorems for space-like graphs with parallel mean curvature and
controlled growth\endtitle
\author Yuxin Dong \endauthor
\thanks {Supported by Zhongdian grant of NSFC}
\endthanks
\abstract {In this paper, we obtain an Ecker-Huisken type result
for entire space-like graphs with parallel mean curvature.}
\endabstract
, {58E20}
\endsubjclass
\endtopmatter
\document
\heading{\bf 1. Introduction}
\endheading
\vskip 0.3 true cm

In 1914, Bernstein proved that the only entire minimal graph in
$R^3$ is a plane. This result was generalized to $R^{m+1}$ for
$m\leq 7$, and higher dimensions and codimensions under various
growth conditions, see [EH], [SXW], [Wa3] and their references. In
1965, Chern [Ch] showed that the only entire graphic hypersurface
in $R^{m+1}$ with constant mean curvature must be minimal.
Therefore we have the corresponding Bernstein type results for
constant mean curvature hypersurfaces. Bernstein type results for
submanifolds in $R^{m+n}$ with parallel mean curvature were also
obtained by some authors (cf. [HJW], [JX1] and [Do]).

In 1968, Calabi [Ca] raised a similar problem for extremal
hypersurfaces in Lorentz -Minkowski space $R_1^{m+1}$ and he
proved that the Bernstein result is true for $2\leq m\leq 4$.
Later, Cheng and Yau [CY] extended Calabi's result to all $m$ as
follows: The only complete extremal space-like hypersurfaces in
$R_1^{m+1}$ are space-like hyperplanes. Recently, Jost and Xin
[JX2] generalized this result to higher codimensional case.

On the other hand, it is important to investigate space-like
constant mean curvature hypersurfaces in $R_1^{m+1}$, which have
interest in relative theory (cf. [MT]). In [Tr], Treibergs showed
that there are many entire space-like graphs with constant mean
curvature besides hyperboloids. Thus Chern type result is no
longer true in this case. It is known that the Gauss map of a
constant mean curvature space-like hypersurface $M$ is a harmonic
map to hyperbolic space. Xin [Xi1] got a Bernstein result by
assuming the boundedness of the Gauss map. Later, [XY] and [CSZ]
extended this result by proving that $M$ must be a space-like
hyperplane if its Gauss image lies in a horoball in the hyperbolic
space. Another natural generalization is to consider a space-like
submanifold in pseudo-Euclidean space $R_n^{m+n}$ with parallel
mean curvature. In [Xi2] the author extended the previous
mentioned result in [Xi1] to higher codimensional case under the
same boundedness assumption on Gauss map.

In this paper, we consider a space-like graphic submanifold
$M=\{(x,f(x)):x\in R^m\}$ in $R_n^{m+n}$ with parallel mean
curvature. Since $M$ is space-like, the induce metric
$(g_{ij})=(\delta _{ij}-\sum_{s=1}^nf_{x_i}^sf_{x_j}^s)$ is
positive definite. Set
$$*\Omega =\{ \sqrt{\det
(I-\sum_{s=1}^nf_{x_i}^sf_{x_j}^s)}\}^{-1}
$$
Our main result is the following:

{\bf Theorem} Let $M^m=(x,f(x))$ be an entire space-like graph in
$ R_n^{m+n}$ with parallel mean curvataure. If the function
$*\Omega $ has growth
$$ \ast \Omega =o(r)\quad \text{as }r\rightarrow
\infty
$$
where $r=\sqrt{\sum_{i=1}^mx_i^2}$, then $M$ is a space-like
$m-$plane.

Our strategy is to establish a Chern-type result for an entire
space-like graph with parallel mean curvature under the growth
condition of $*\Omega $. Then the result follows immediately from
[CY] and [JX2]. Notice that if $n=1$ , we have
$$
\ast \Omega =\frac 1{\sqrt{1-|\nabla f|^2}}
$$
Therefore the growth condition of $*\Omega $ is similar to that
one given by Ecker-Husken [EH] for minimal graphic hypersurfaces
in $R^{m+1}$. The above result may be regarded as an Ecker-Huisken
type result for space-like graphs with parallel mean curvature. By
calculating the quantity $*\Omega $ of the hyperboloid, we will
see that the growth condition is optimal. In [Do], the author uses
a similar method to establish some Bernstein type results for
submanifolds in Euclidean space with parallel mean curvature.

\heading{\bf 2. Preliminaries}
\endheading
\vskip 0.3 true cm

In this section, we will generalize Chern's method [Ch] to our setting. Let $%
R_n^{m+n}$ be an $(m+n)-$dimensional pseudo-Euclidean space of
index $n$, namely the vector space $R^{m+n}$ endowed with the
metric
$$
(\,\,,\,)=(dx_1)^2+\cdots +(dx_m)^2-(dx_{m+1})^2-\cdots
-(dx_{m+n})^2\tag{1}
$$
The standard Euclidean metric of $R^{m+n}$ will be denoted by
$(\,,\,)_E$. For a vector $v$ in $R^{m+n}$, we will use the
notations $|v|$ and $|v|_E$ to denote the norms of $v$ with
respect to $(\,,\,)$ and $(\,,\,)_E$ respectively.

Let $z:M^m\rightarrow R_n^{m+n}$ be a space-like immersion of an
oriented $ m- $dimensional manifold into $R_n^{m+n}$. We will
regard $z$ as a vector-valued function on $M$. Choose a local
Lorentzian frame field $ \left\{
e_1,...,e_m,e_{m+1},...,e_{m+n}\right\}$ such that $
\{e_{m+1},...,e_{m+n}\}$ is a normal frame field of $M$. Throught
this paper, we agree with the following indices:
$$
\aligned &1\leq A,B,C...\leq m+n\\
 1\leq i,j,k,...&\leq m,\quad m+1\leq \alpha ,\beta ,\gamma ,....\leq m+n
 \endaligned \tag{2}
$$
Write
$$
\aligned
&dz=\sum_A\omega _Ae_A \\
&de_A=\sum_B\omega _{AB}e_B\endaligned\tag{3}
$$
Therefore $\{\omega _i\}$ is a dual frame field of $\left\{
e_i\right\} $ and $\omega _\alpha =0$ on $M$. The induced
Riemannian metric of $M$ is then given by $ds_M^2=\sum_i\omega
_i^2$. By Cartan's lemma, we have
$$
\omega _{\alpha i}=\sum_kh_{\alpha ij}\omega _j,\quad h_{\alpha
ij}=h_{\alpha ji} \tag{4}
$$
where $h_{\alpha ij}$ are components of the second fundamental
form of $M$ in $R_n^{m+n}$. The mean curvature vector of $M$ is
defined by
$$
\vec{H}=\frac 1m\sum_{\alpha,k}h_{\alpha kk}e_\alpha \tag{5}
$$
If $\nabla ^{\bot }\vec{H}=0$, $M$ is said to have parallel mean
curvature. If $\vec{H}=0$, $M$ is called an extremal spacelike
submanifold.

Now we consider a space-like graph $M=\{$ $(x,f(x)):x\in D\subset
R^m\}$ in $R_n^{m+n}$ with parallel mean curvture $\vec{H}$, where
$D$ is a compact domain with smooth boundary $\partial D$.
Obviously $H=\sqrt{-(\vec{H},\vec{H})}$ is a nonnegative constant.

Let $\Omega =dx^1\wedge \cdots \wedge dx^m$ be the parallel
$m-$form on $R_n^{m+n}$ and let $\{a_1,...,a_{m+n}\}$ be an
oriented Lorentzian basis of $R_n^{m+n}$ such that
$\{a_i\}_{i=1}^m$ is an oriented orthonormal basis of $R^m$ . If
$H>0$, we have a global future-directed normal vector field
$e_H=H^{-1}\vec{H}$. Therefore we may define a global $m-$ form on
$M$ as follows:
$$
\Phi =(m-1)!\sum_{i=1}^md(a_1,z)\wedge \cdots \wedge
d(a_{i-1},z)\wedge d(a_i,e_H)\wedge d(a_{i+1},z)\wedge \cdots
\wedge d(a_m,z)\tag{6}
$$
Clearly $\Phi$ is independent of the choice of the oriented
orthogonal basis $\left\{ a_i\right\} _{i=1}^m$ in $R^m$.

For any $p\in M$ then the differential of $f$ is a linear map from
$R^m$ to $R^n$. As in [Wa1], we can use singular value
decomposition to find orthonormal bases $\left\{ a_i\right\}
_{i=1}^m$ for $R^m$ and $\left\{ a_\alpha \right\} _{\alpha
=m+1}^{m+n}$ for $R^m$ such that
$$
df(a_i)=\lambda _ia_{m+i}\tag{7}
$$
for $i=1,...,m$. Notice that $\lambda _i=0$ if $i>\min \left\{
m,n\right\}$. Then we have
$$
(a_i,a_j)=\delta _{ij},\ (a_i,a_\alpha )=0,\ (a_\alpha ,a_\beta
)=-\delta _{\alpha \beta }
$$
Therefore we have a Lorenztian basis $\left\{ e_A\right\} $ at $p$
given by
$$
\left\{ e_i=\frac 1{\sqrt{1-\lambda _i^2}}(a_i+\lambda
_ia_{m+i})\right\} _{i=1,...,m}\in T_pM  \tag{8}
$$
and
$$
\left\{ e_\alpha =\frac 1{\sqrt{1-\lambda _{\alpha
-m}^2}}(a_\alpha +\lambda _{\alpha -m}a_{\alpha -m})\right\}
_{\alpha =m+1}^{m+n}\in T_p^{\bot }M \tag{9}
$$
By definition $*\Omega =\Omega (e_1,...,e_m)$, and thus we have
$$
\ast \Omega =\frac 1{\sqrt{\Pi _{i=1}^m(1-\lambda _i^2)}}\tag{10}
$$
\proclaim{Lemma 1} Under the above notations, we have
$$
\Phi =m!H(*\Omega )\omega ^1\wedge \cdots \wedge \omega ^m
$$
where $\omega ^1\wedge \cdots \wedge \omega ^m$ is volume form of
$M$.
\endproclaim
\demo{Proof} Using(3), (8) and (9), we have from (6) the
following:
$$
\aligned \Phi=&(m-1)!\sum_{i=1}^m(a_1,dz)\wedge \cdots \wedge
(a_{i-1},dz)\wedge(a_i,de_H)\wedge (a_{i+1},dz)\wedge \cdots \wedge (a_m,dz)\\
=&(m-1)!\sum_{i=1}^m\frac 1{\sqrt{1-\lambda _1^2}}\omega
^1\wedge\cdots \wedge \frac 1{\sqrt{1-\lambda _{i-1}^2}}\omega
^{i-1}\wedge \left( \frac{ -h_{ii}^H}{\sqrt{1-\lambda _i^2}}\omega
^i\right)\wedge\\
&\frac 1{\sqrt{1-\lambda _{i+1}^2}}\omega
^{i+1}\wedge\cdots \wedge \frac 1{\sqrt{1-\lambda _m^2}}\omega ^m \\
=&-(m-1)!(*\Omega )(\sum_ih_{ii}^H)\omega ^1\wedge \cdots \wedge
\omega ^m\\
=&m!(*\Omega )H\omega ^1\wedge \cdots \wedge \omega ^m
\endaligned
$$
where $\sum_ih_{ii}^H=<\sum_\alpha h_{\alpha ii}e_\alpha
,e_H>=-mH$. This proves the Lemma. \qed
\enddemo

We may write
$$
\Phi =(m-1)!d\alpha  \tag{11}
$$
where
$$
\left. \alpha =\sum_i(-1)^{i-1}(a_i,e_H)d(a_1,z)\wedge \cdots
\wedge d(a_{i-1},z)\wedge d(a_{i+1},z)\wedge \cdots \wedge
d(a_m,z)\right.\tag{12}
$$
Applying the Stokes Theorem to (11), we get
$$
\left. mH\int_M(*\Omega )\omega ^1\wedge \cdots \wedge \omega
^m=\int_{\partial M}\alpha .\right.  \tag{13}
$$

We project $z(M)$ orthogonally into the $m-$plane spaned by
$\left\{ a_i\right\} _{i=1}^m$. If $z^{\prime }(p)$ is the image
point of $z(p)$, $ p\in M$, under this orthogonal projection, we
have
$$
z^{\prime }=z+\sum_{\alpha =m+1}^{m+n}(a_\alpha ,z)a_\alpha
\tag{14}
$$
Let $\Psi $ be a non-zero differential form on $\partial
M=\{(x,f(x)):x\in
\partial D\}$, defined locally. Using this form, the elements of volume of $
z^{\prime }(\partial M)$, $z(\partial M)$ may be expressed
respectively as $P\Psi $, $Q\Psi $ with $P\geq 0$ and $Q\geq 0$.
We write
$$
\omega _{i_1}\wedge \cdots \wedge \omega
_{i_{m-1}}=p_{i_1,...,i_{m-1}}\Psi \tag{15}
$$
on $\partial M$.

By a direct computation, we have
$$
\aligned \frac 1{(m-1)!}\underset{m-1}\to {\underbrace{dz\wedge
\cdots \wedge dz}}&=\frac 1{(m-1)!}(\sum \omega
_{i_1}e_{i_1})\wedge \cdots \wedge (\sum \omega
_{i_{m-1}}e_{i_{m-1}}) \\
&=\sum_ip_{1,...,i-1,i+1,...,m}(e_1\wedge \cdots \wedge
e_{i-1}\wedge e_{i+1}\wedge \cdots \wedge e_m)\Psi
\endaligned
$$
so that
$$
Q^2=\sum_ip_{1,...,i-1,i+1,...,m}^2.  \tag{16}
$$
Using (8), (9) and (10), we get:
$$
\aligned \alpha &=\sum_{i=1}^m(-1)^{i-1}(a_i,e_H)(a_1,dz)\wedge
\cdots \wedge
(a_{i-1},dz)\wedge (a_{i+1},dz)\wedge \cdots \wedge (a_m,dz) \\
&=\sum_{i=1}^m(-1)^{i-1}(a_i,e_H)\prod_{j\not =i}\frac
1{\sqrt{1-\lambda _j^2}}\omega ^1\wedge \cdots \wedge \omega
^{i-1}\wedge \omega ^{i+1}\wedge
\cdots \wedge \omega ^m \\
&=\sum_{i=1}^m(-1)^{i-1}(a_i,e_H)p_{1...i-1,i+1,...,m}\Psi \\
&=(*\Omega )\sum_{i=1}^m(-1)^i\lambda _i\xi
_{m+i}p_{1...i-1,i+1,...,m}\Psi \endaligned  \tag{17}
$$
where $\xi _{m+k}=<e_H,e_{m+k}>$ if $k\leq \min \left\{
m,n\right\} $ and $\xi _{m+k}=0$ if $k>\min \left\{ m,n\right\}$.
Obviously $\sum_{i=1}^m\xi _{m+i}^2\leq 1$. Since $|df|<1$, we get
from (17) and the Cauchy-Schwarz inequality that
$$
|R|\leq (*\Omega )Q \tag{18}
$$

Next, we will show that if $M$ is a space-like hypersurface, there
is a nice formula relating the quantities $P,Q$ and $R$. When
$n=1$, (12) is simplified to
$$
\alpha =vp_{2,...,m}  \tag{19}
$$
where $v=\lambda _1/\sqrt{1-\lambda _1^2}=\lambda _1(*\Omega )$;
and thus
$$
R^2=v^2p_{2,...,m}^2 \tag{20}
$$
Write $a=a_{m+1}$. Then (14) becomes
$$
x^{\prime }=x+(a,x)a  \tag{21}
$$
From (9) and (10), we easily derive
$$
(a,e_i)=-\delta _{i1}v,\quad (a,e_{m+1})=-*\Omega\tag{22}
$$
and
$$
a=-ve_1+*\Omega e_{m+1} \tag{23}
$$
To determine $P$, we compute $\frac 1{(m-1)!} \underset{m-1}\to
{\underbrace{ dz^{\prime }\wedge \cdots \wedge dz^{\prime }}}$ as
follows:
$$
\aligned &\frac 1{(m-1)!}\underset{m-1}\to {\underbrace{dz^{\prime
}\wedge\cdots\wedge dz^{\prime}}}\\
&=\frac 1{(m-1)!}\{\sum_{i_1}e_{i_1}\omega _{i_1}-v\omega
_1a\}\wedge\cdots\wedge\{\sum_{i_{m-1}}e_{i_{m-1}}\omega_{i_{m-1}}-v\omega_1a\}\\
&=\frac 1{(m-1)!}\sum\omega _{i_1}\wedge\cdots\omega_{i_{m-1}}
(e_{i_1}\wedge\cdots\wedge e_{i_{m-1}})-\frac v{(m-1)!}
\underset{1\leq s\leq m-1}\to {\sum }\omega_{i_1}\wedge\cdots\wedge\\
&\omega _{i_{s-1}}\wedge\omega _1\wedge \omega _{i_{s+1}}\wedge
\cdots \wedge\omega_{i_{m-1}}(e_{i_1}\wedge\cdots\wedge a\wedge \cdots\wedge e_{i_{m-1}})\\
&=\frac 1{(m-1)!}\sum \omega _{i_1}\wedge\cdots \omega_{i_{m-1}}
(e_{i_1}\wedge\cdots\wedge e_{i_{m-1}})+\frac{v^2}{(m-1)!}
\underset{1\leq s\leq m-1}\to {\sum}(-1)^{s-1}\omega _1\wedge \\
&\omega _{i_1}\wedge\cdots\wedge \omega _{i_{s-1}}\wedge
\omega_{i_{s+1}}\cdots \wedge \omega _{i_{m-1}} (e_{i_1}\wedge
\cdots \wedge e_1\wedge \cdots \wedge e_{i_{m-1}})\\
&-\frac {v*\Omega}{(m-1)!}\underset{1\leq s\leq
m-1}\to{\sum}(-1)^{s-1}\omega _1\wedge\omega _{i_1}\wedge \cdots
\wedge \omega _{i_{s-1}}\wedge \omega_{i_{s+1}}\cdots\wedge \omega
_{i_{m-1}}(e_{i_1}\wedge\cdots \wedge\\
&e_{m+1}\wedge \cdots \wedge e_{i_{m-1}})
\endaligned
$$
So the coefficient of $\Psi $ in $\frac 1{(m-1)!}\underset{m-1}\to
{\underbrace{dz^{\prime }\wedge \cdots \wedge dz^{\prime }}}$ is
$$
\aligned &\underset{i_1<...<i_{m-1}}\to {\sum}
p_{i_1....i_{m-1}}e_{i_1}\wedge\cdots\wedge
e_{i_{m-1}}+v^2\underset{1<i_2<...<i_{m-1}}\to {\sum}
p_{1i_2...i_{m-1}}e_1\wedge e_{i_2}\wedge \cdots \wedge e_{i_{m-1}}\\
&-(-1)^mv(*\Omega)\underset{i_2<...<i_{m-1}}\to {\sum}
p_{1i_2...i_{m-1}}e_{i_2}\wedge \cdots \wedge e_{i_{m-1}}\wedge
e_{m+1}
\endaligned
$$
It follows that
$$
\aligned P^2 &=\underset{1<i_1<...<i_{m-1}}\to {\sum
}p_{i_1....i_{m-1}}^2+(1+v^2)^2
\underset{1<i_2<...<i_{m-1}}\to {\sum }p_{1i_2...i_{m-1}}^2 \\
&-v^2(*\Omega )^2\underset{1<i_2<...<i_{m-1}}\to {\sum }p_{1i_2...i_{m-1}}^2 \\
&=\underset{i_1<...<i_{m-1}}\to {\sum
}p_{i_1....i_{m-1}}^2+v^2\underset{1<i_2<...<i_{m-1}}\to{\sum }p_{1i_2...i_{m-1}}^2 \\
&=\underset{i_1<...<i_{m-1}}\to{\sum}p_{i_1....i_{m-1}}^2+v^2(\underset{
i_1<...<i_{m-1}}\to {\sum }p_{i_1....i_{m-1}}^2-p_{2...m}^2)\\
&=(1+v^2)Q^2-R^2
\endaligned
$$
since $(e_{m+1},e_{m+1})=-1$ and $v^2-(*\Omega )^2=-1$. Thus we
have
$$
P^2+R^2=(*\Omega )^2Q^2 \tag{24}
$$

\heading{\bf 3. Bernstein-Type Theorems}
\endheading
\vskip 0.3 true cm

In this section, we take $D=\left\{ x\in R^m:\sum_{i=1}^mx_i^2\leq
r\right\} $. As before, let
$$
M=\left\{ (x,f(x):x\in D\right\}
$$
be a space-like graph in $R_n^{m+n}$ with paralle mean curvature.
\proclaim {Lemma 1} On $\partial M$, we have
$$
Q\leq P \tag{25}
$$
In particular, if $n=1$, i.e., $M$ is a space-like hypersurface,
then we have
$$
Q\leq \sqrt{(*\Omega )^{-2}+|df(\eta _m)|^2}P\tag{26}
$$
where $|df(\eta _m)|$ may be regarded as the radial singular value
of the map $f$.
\endproclaim
\demo{Proof} Choose an orthonormal basis $\left\{ \eta _1,...,\eta
_m\right\} $ at $q\in \partial D$ in $R^m$ such that $\eta _m$ is
a normal vector of $\partial D$. We have the corresponding tangent
vectors of the graph $M$ at $(q,f(a))$
$$
\xi _{i\text{ }}=(\eta _i,df(\eta _i)),\quad i=1,...,m
$$
It is easy to see
$$ |\xi _1\wedge \cdots \wedge \xi _m|=(*\Omega
)^{-1}  \tag{27}
$$
and
$$
|\xi _1\wedge \cdots \wedge \xi _{m-1}|=Q/P  \tag{28}
$$
Write
$$
\xi _i=\widetilde{\eta _i}+\widetilde{df(\eta _i)}
$$
where $\widetilde{\eta _i}=(\eta _i,0)$ and $\widetilde{df(\eta
_i)} =(0,df(\eta _i))$. Therefore
$$
|\xi _i|^2=|\widetilde{\eta _i}|_E^2-|\widetilde{df(\eta _i)}|_E^2=1-|%
\widetilde{df(\eta _i)}|_E^2\leq 1
$$
and thus
$$
|\xi _1\wedge \cdots \wedge \xi _{m-1}|\leq 1\tag{29}
$$
From (28) and (29), we have (25).

Now assume that $n=1$. Obviously $\xi _i$, $i=1,...,m-1$, and
$\widetilde{f(\eta _m)}$ are tangent to the cylinder
$$
C_D=\{(x_1,....,x_{m+n})\in R^{m+n}:\sum_{i=1}^mx_i^2=r\}
$$
and the horizontal vector $\widetilde{\eta _m}$ is orthogonal to
$C_D$ at the point $(q,f(q))$. It follows that
$$
\aligned
(*\Omega )^{-2}&=|\xi _1\wedge \cdots \wedge \xi _m|^2 \\
&=|\xi _1\wedge \cdots \wedge \xi _{m-1}\wedge \widetilde{\eta
_m}+\xi
_1\wedge \cdots \wedge \xi _{m-1}\wedge \widetilde{df(\eta _m)}|^2 \\
&=|\xi _1\wedge \cdots \wedge \xi _{m-1}\wedge \widetilde{\eta
_m}|^2+|\xi
_1\wedge \cdots \wedge \xi _{m-1}\wedge \widetilde{df(\eta _m)}|^2 \\
&=|\xi _1\wedge \cdots \wedge \xi _{m-1}\wedge \widetilde{\eta
_m}|^2+| \widetilde{\eta _1}\wedge \cdots \wedge \widetilde{\eta
_{m-1}}\wedge
\widetilde{df(\eta _m)}|^2 \\
&\geq |\xi _1\wedge \cdots \wedge \xi
_{m-1}|^2-|\widetilde{df(\eta _m)}|_E^2
\endaligned
$$
i.e.,
$$
Q^2/P^2=|\xi _1\wedge \cdots \wedge \xi _{m-1}|^2\leq (*\Omega
)^{-2}+|\widetilde{df(\eta _m)}|_E^2
$$
This gives (26).\qed
\enddemo

We recall the following \proclaim{Theorem A}$([CY], [JX])$ Let $M$
be an extremal space-like $m-$ submanifold in $R_n^{m+n}$. If $M$
is closed with respect to the Euclidean topology, then $M$ has to
be a space-like $m-$plane.
\endproclaim

\remark {Remark} Obviously, an entrie graph is closed with respect
to the Euclidean topology. Hence we know that any entire extremal
space-like graph must be a space-like $m-$plane.
\endremark

\proclaim{Theorem 1}Let $M^m=(x,f(x))$ be an entire space-like
graph in $R_n^{m+n}$ with parallel mean curvature. If $*\Omega $
has the following growth
$$
\ast \Omega =o(r)\quad \text{as }r\rightarrow \infty \tag{30}
$$
where $r=\sqrt{\sum_{i=1}^mx_i^2}$, then $M$ is a space-like
$m-$plane.
\endproclaim

\demo {Proof} Let $M_r=\left\{ (x,f(x)):x\in D_r\subseteq
R^m\right\} $, where $D_r$ denotes the closed ball of radius $r$
centered at the origin in $R^m$. From (13), (18) and Lemma 1, we
have
$$
\aligned mH\int_{M_r}(*\Omega )\omega ^1\wedge \cdots \wedge
\omega ^m&\leq\int_{\partial M_r}R\Psi \\
&\leq \int_{\partial M_r}*\Omega P\Psi \\
&\leq \sup_{\partial D_r}\{*\Omega \}Vol(\partial D_r)
\endaligned
$$
i.e.,
$$
mHVol(D_r)\leq \sup_{\partial D_r}\{*\Omega \}Vol(\partial D_r)
$$
Thus
$$
H\leq C\frac{\sup_{\partial D_r}\{*\Omega \}}r
$$
where $C$ is a universal constant. Let $r\rightarrow \infty $. It
follows that $H\equiv 0$. Hence we may complete the proof by
Theorem A.\qed
\enddemo

For space-like hypersurfaces, we may give a more delicate growth
condition to ensure the above result.

\proclaim{Proposition 2}Let $M^m=(x,f(x))$ be an entire space-like
hypersurface in $R_1^{m+1}$ with constant mean curvature. If
$$
\sup_{\partial D_r}\{|df(\eta _m)|_E*\Omega \}=o(r)
$$
where $r=\sqrt{\sum_{i=1}^mx_i^2}$, then $M$ is a space-like
$m-$plane.
\endproclaim

\demo{Proof} From (13), (25) and Lemma 1, we have
$$
\aligned mH\int_{M_r}(*\Omega )\omega ^1\wedge \cdots \wedge
\omega ^m&\leq
\int_{\partial M_r}R\Psi \\
&\leq \int_{\partial M_r}\sqrt{(*\Omega )^2Q^2-P^2}\Psi \\
&\leq \int_{\partial M_r}|df(\eta _m)|_E(*\Omega )P\Psi \\
&\leq \sup_{\partial D_r}\{|df(\eta _m)|_E*\Omega \}Vol(\partial
D_r) \endaligned $$ By the same argument as in Theorem 1, we prove
the proposition. \qed
\enddemo

Let's consider a typical example of space-like graphs in
$R_1^{m+1}$with constant mean curvature.

\example{Example 1}The hyperboloid is defined by
$$
\aligned H_{-1}^m &=\{(x_1,...,x_n,x_{n+1})\in
R_1^{m+1}:\sum_{i=1}^mx_i^2-x_{m+1}^2=-1,x_{m+1}\geq 0\} \\
&=\{(x,f(x)):f=\sqrt{1+\sum_{i=1}^mx_i^2},x\in R^m\}
\endaligned
$$
By a direct computation, we have
$$
\ast \Omega =\frac 1{\sqrt{1-|\nabla
f|^2}}=\sqrt{1+\sum_{i=1}^mx_i^2} =O(r)\tag{31}
$$
From (31), we see that the growth condition in Theorem 1 is
optimal.
\endexample

\proclaim{Theorem B}$([XY], [CSZ])$ Let $M^m=(x,f(x))$ be a
complete space-like hypersurface in $R_1^{m+1}$ with constant mean
curvature. If the image of the Gauss map $\gamma :M\rightarrow
H^m(-1)$ lies in a horoball in $H^m(-1)$, then $M$ must be a
space-like hyperplane.
\endproclaim

It is known that every complete spacelike hypersurface in
$R_1^{m+1}$ is spatially entire (cf. [AM]). To compare Theorem 1
with Theorem B, we hope to find the equivalent restriction on the
function $*\Omega $, if the image of $\gamma $ lies in a horoball.

Let $M=(x,f(x))$ be a space-like graphic hypersurface in
$R_1^{m+1}$. Its Gauss map $\gamma $ is given by
$$
\aligned
\gamma :M&\longrightarrow H_{-1}^m \\
x&\longmapsto \frac 1{\sqrt{1-|\nabla
f|^2}}(f_{x_1},...,f_{x_m},1)=*\Omega (f_{x_1},...,f_{x_m},1)
\endaligned \tag{32}
$$
where $H_{-1}^m$ is the hyperboloid endowed with the induced
metric from $ R_1^{m+1}$. Obviously the Gauss image of $M$ is
bounded in $H_{-1}^m$ if and only if $*\Omega $ is bounded. This
also holds true for higher codimensional case(cf. [Xi2]).

It is easier to use the upper half-space model $H^m$ of the
hyperbobolic space for describing horoballs. We consider the
following maps
$$
\aligned
h_1:H_{-1}^m&\longrightarrow B^m \\
(x_1,....,x_m,x_{m+1})&\longmapsto
(\frac{x_1}{1+x_{m+1}},....,\frac{x_m}{1+x_{m+1}})
\endaligned\tag{33}
$$
and
$$
\aligned
h_2:B^m&\longrightarrow H^m=\left\{ (y_1,...,y_m)\in R^m:y_m>0\right\} \\
p&\longmapsto 2\frac{p-p_0}{|p-p_0|^2}-(0,...,0,1)
\endaligned \tag{34}
$$
where $p_0=(0,...,-1)$ and $H^m$ is endowed with the metric $
g=y_m^{-2}(dy_1^2+\cdots +dy_m^2)$. The set $\left\{
(y_1,...,y_m)\in H^m:y_m>c>0\right\}$ for any positive constant
$c$ is a horoball in $H^m$. It is known that $h_2\circ
h_1:H_{-1}^m\rightarrow H^m$ is an isomorphism. From (32), (33)
and (34), we may get the $m-$th component of $h_2\circ h_1\circ
\gamma $ as follows:
$$
(h_2\circ h_1\circ \gamma )_m=\frac 1 {(1+f_{x_m})*\Omega }
$$
So the condition $y_m>c>0$ is equivalent to
$$
(1+f_{x_m})*\Omega < \frac 1 c \tag{35}
$$
Note that $f_{x_m}$ may be replaced by any $f_{x_i}$ or $v(f)$
which denotes the derivative in any fixed unit direction $v$ in
$R^m$. Obviously, if there exists a sequence of points $\left\{
p_k\right\}$ such that $*\Omega (p_k)\rightarrow \infty $, then
$(f_{x_m})(p_k)\rightarrow -1$. Therefore (35) implies that all
`bad singular directions' approach one direction, i.e., $\partial
/\partial x_m$.

Since $*\Omega =(\sqrt{1-|\nabla f|^2})^{-1}$, we see that the
growth condition in Theorem 1 is very much like that one given by
Ecker-Huisken in [EH] for a minimal graphic hypersurface in the
Euclidean space $R^{m+1}$. Hence Theorem 1 may be regarded as an
Ecker-Huisken type result.

\vskip 0.3 true cm {\bf Acknowledgments:} The author would like to
thank Professors C.H. Gu and H.S. Hu for their constant
encouragement and helpful comments.

\vskip 1 true cm 
\vskip 1 true cm

Institute of Mathematics

Fudan University, Shanghai 200433

P.R. China

And

Key Laboratory of Mathematics

for Nonlinear Sciences

Ministry of Education

\vskip 0.2 true cm
yxdong\@fudan.edu.cn

\enddocument